\documentclass{amsart}
\usepackage{amsmath,amsthm,amssymb,amscd,amsfonts}
\usepackage{enumerate}
\makeatletter
\@namedef{subjclassname@2010}{%
  \textup{2010} Mathematics Subject Classification}
\makeatother
\theoremstyle{definition}

\newtheorem{theo}{Theorem}[section]
\newtheorem{cor}[theo]{Corollary}
\newtheorem{lemm}[theo]{Lemma}
\newtheorem{prop}[theo]{Proposition}
\theoremstyle{definition}
\newtheorem{defi}[theo]{Definition}
\theoremstyle{definition}
\newtheorem{ex}[theo]{Example}
\newtheorem{rem}[theo]{Remark}
\numberwithin{equation}{section}

\title{ Weakly Connes amenable dual Banach algebras }

\author{Amin Mahmoodi}
\address{Department of Mathematics, Central Tehran Branch, Islamic Azad University, Tehran, Iran, e-mail: {\tt a\_mahmoodi@iauctb.ac.ir}}

\thanks{} 

\begin{document}
\maketitle

\setcounter{section}{0}
\begin{abstract}
We shall develop a notion of amenability for dual Banach algebras, namely
weak Connes amenability, which will play the role that weak amenability does for usual Banach algebras.

\end{abstract}

\section{Introduction}

Amenable Banach algebras were introduced by B. E. Johnson in
\cite{joh}. There are several variants of amenability, two of the
most notable are weak amenability and Connes amenability. The
concept of weak amenability for Banach algebras was introduced by B.
E. Johnson in \cite{joh1}, it generalizes that introduced by W. G.
Bade, P. C. Curtis, and H. G. Dales for commutative Banach algebras
in \cite{bcd}. The notion of Connes amenability systematically
introduced by V. Runde in \cite{r1}. We recall the definitions in
Definitions \ref{1.1} and \ref{1.2} below.

The purpose of this paper is to study a new notion of amenability for dual Banach algebras. The organization of the paper is as
 follows. In Section 2, we recall some background notations and definitions.

  In section 3, weak Connes amenability for dual Banach algebras is introduced and some basic and hereditary properties are given.  It is shown that the corresponding class of such algebras includes all Connes amenable dual Banach algebras (Theorem \ref{0.2}), as well as all
weakly amenable dual Banach algebras  (Theorem \ref{2.12}). It is proved that commutative, pseudo Connes amenable dual Banach
algebras are always weakly Connes amenable (Corollary \ref{2.8}).  We study weak Connes amenability of direct sums of dual Banach
algebras (Theorem \ref{5.2}). Weak Connes amenability of the enveloping dual
Banach algebras is also discussed (Corollary \ref{2.11}).

In section 4 we verify this new notion for some certain algebras. Examples are given to
distinguish between the new notion and the classical concepts of amenability. In
particular, we present some weakly Connes amenable dual Banach
algebras which are neither Connes amenable nor weakly
amenable (Theorems \ref{5.4} and \ref{5.3}).

The author thanks an anonymous referee for suggestions, in response to a previous version of this work, which improved the paper.

\section{Preliminaries}

Suppose that $\mathfrak{A}$ is a Banach algebra. It is known that the projective tensor product
$\mathfrak{A} \hat{\otimes} \mathfrak{A}$ is a Banach
$\mathfrak{A}$-bimodule in the canonical way. There is a continuous
linear $\mathfrak{A}$-bimodule homomorphism $ \pi : \mathfrak{A}
\hat{\otimes} \mathfrak{A} \longrightarrow \mathfrak{A}$ such that $
\pi ( a \otimes b) = ab$ for $a,b \in \mathfrak{A}$. If $E$ is a
Banach $\mathfrak{A}$-bimodule, then so is the dual space $E^*$.  A continuous linear map
$D:\mathfrak{A} \longrightarrow E$ is a \textit{derivation} if it
satisfies $ D(ab) = D(a) \ . \ b + a \ . \ D(b) $ for all $a,b \in
\mathfrak{A}$. We call $D$ \textit{inner} if there is $x \in
E$ such that $ D(a) = ad_x(a) := a \ . \ x - x \ . \ a$ for every $a \in
\mathfrak{A}$.

\begin{defi} \label{1.1} A Banach algebra $\mathfrak{A}$ is \textit{weakly
amenable} if every derivation $D: \mathfrak{A} \longrightarrow
\mathfrak{A}^*$ is inner.
\end{defi}

Let $X$ be a Banach space. We simply denote by $wk$ and $w^*$, the $\sigma( X, X^*)$-topology and the $\sigma( X^*, X)$-topology on $X$ and $X^*$, respectively.

Let $\mathfrak{A}$ be a Banach algebra. A Banach
$\mathfrak{A}$-bimodule $E$ is \textit{dual} if there is a closed
submodule $E_*$ of $E^*$ such that $E = (E_*)^*$. We call $E_*$ the
\textit{predual} of $E$. A Banach algebra $\mathfrak{A}=(\mathfrak{A}_*) ^*$ is
\textit{dual} if it is dual as a Banach $\mathfrak{A}$-bimodule. Equivalently, a Banach algebra $\mathfrak{A}$ is dual if it is a dual Banach space such that its multiplication is separately continuous in the $w^*$-topology.

 Let
$\mathfrak{A}$ be a dual Banach algebra, and let $E$ be a dual
Banach $\mathfrak{A}$-bimodule. Then we say $E$ is \textit{normal}
if the module actions of $\mathfrak{A}$ on $E$ are $w^*$-$w^*$
continuous.
\begin{defi} \label{1.2} A dual Banach algebra  $\mathfrak{A}$ is \textit{Connes amenable}
if for every normal, dual Banach $\mathfrak{A}$-bimodule $E$, every
$w^*$-$w^*$ continuous derivation $D : \mathfrak{A} \longrightarrow
E$ is inner.
\end{defi}
 Let $\mathfrak{A} = (\mathfrak{A}_*)^*$ be a
dual Banach algebra and let $E$ be a Banach $\mathfrak{A}$-bimodule.
We write $\sigma wc(E)$ for the set of all elements $ x \in E$ such
that the maps $$ \mathfrak{A} \longrightarrow E \ \ , \ \ \ a
\longmapsto \left \{
\begin{array}{ll}
                        a \ . \ x \\
                        x \ . \ a
                      \end{array} \right. \ , $$ are $w^*$-$wk$
continuous. It is well known that $\sigma wc(E)$ is a closed submodule of $E$, and  $\sigma wc(E)^*$ is always normal.
It is shown in \cite [Corollary 4.6]{r3}, that $
\pi^*(\mathfrak{A}_*) \subseteq \sigma wc(\mathfrak{A} \hat{\otimes}
\mathfrak{A} )^*$ . Taking adjoints, we can extend $\pi$ to an
$\mathfrak{A}$-bimodule homomorphism $\pi_{\sigma wc}$ from $ \sigma
wc((\mathfrak{A} \hat{\otimes} \mathfrak{A})^*)^*$ to
$\mathfrak{A}$. A $\sigma wc$-\textit{virtual diagonal} for a dual
Banach algebra $\mathfrak{A}$ is an element $M \in \sigma
wc((\mathfrak{A} \hat{\otimes} \mathfrak{A})^*)^*$ such that $ a \ .
\ M = M \ . \ a$ and $ a \pi_{\sigma wc} (M) = a$ for $a \in
\mathfrak{A}$. It is known that Connes amenability of $\mathfrak{A}$
is equivalent to existence of a $\sigma wc$-virtual diagonal for
$\mathfrak{A}$ \cite {r3}. 

Our comprehensive references on amenability of Banach algebras are \cite{dal, r2}.

\section{Weak Connes amenability}
Let $\mathfrak{A} = (\mathfrak{A}_*)^*$ be a dual Banach algebra, and let $E$ be a Banach
$\mathfrak{A}$-bimodule.  We denote by $j_E : E^*  \longrightarrow
\sigma wc (E)^* $ the adjoint of the inclusion map $  \sigma wc (E) \hookrightarrow E$. It is clear that $j_E$ is an $\mathfrak{A}$-bimodule homomorphism and  $ \langle  x , j_E (f)  \rangle =
f (x)$ for all $x \in \sigma wc (E)$ and $ f \in E^*$.

We start with the following observation.

\begin{lemm} \label{0.1} Let $\mathfrak{A} = (\mathfrak{A}_*)^*$ be a dual Banach algebra, let $E$ be a Banach
$\mathfrak{A}$-bimodule, and let $D :\mathfrak{A} \longrightarrow E^*$ be a derivation. Then:

$(i)$ The map  $j_{E} \circ D : \mathfrak{A} \longrightarrow   \sigma wc (E )^*$ is a derivation;

$(ii)$ If $D$ is inner, then so is  $j_{E} \circ D$. Furthermore,  $j_{E} \circ  ad_f = ad_{j_{E}(f)}$ for each $f \in E^*$;

$(iii)$  For every $f \in E^*$,  $j_{E} \circ  ad_f : \mathfrak{A} \longrightarrow   \sigma wc (E )^*$ is $w^*$-$w^*$ continuous.
\end{lemm}
{\bf Proof.} $(i)$ For every $a, b \in \mathfrak{A} $ and $ x \in  \sigma wc (E )$ we have
\begin{align*}
\langle  x , a \cdot (j_{E} \circ D  ) (b) +  (j_{E} \circ D  ) (a) \cdot b  \rangle &= \langle  x ,  j_{E}( a \cdot D (b) ) + 
 j_{E} ( D (a) \cdot b )  \rangle \\&= \langle  x ,  a \cdot D (b)  +  D (a) \cdot b  \rangle \\&=  \langle  x ,  D (a b)  \rangle =
 \langle  x , ( j_{E} \circ D )  (a b)  \rangle \ ,
\end{align*}
and hence  $j_{E} \circ D$ is a derivation. 

$(ii)$ This is a simple calculation.

$(iii)$ It follows from $(ii)$ and normality of $ \sigma wc (E )^* $.  \qed

\begin{defi} \label{2.1}  A dual Banach algebra $\mathfrak{A} = (\mathfrak{A}_*)^*$ is \textit{weakly Connes amenable} if for every derivation $D
:\mathfrak{A} \longrightarrow \mathfrak{A}^*$ such that  $j_{\mathfrak{A}} \circ D : \mathfrak{A} \longrightarrow   \sigma wc (\mathfrak{A} )^*$ is $w^*$-$w^*$ continuous, derivation  $j_{\mathfrak{A}} \circ D$ is inner.
\end{defi}

\begin{rem} \label{7.1} For a dual Banach algebra $\mathfrak{A} $, $\sigma wc (\mathfrak{A} )$ is a closed (two-sided) ideal of $\mathfrak{A} $. To see this, take $a \in \sigma wc (\mathfrak{A} )$, $b \in \mathfrak{A} $, and let $ c_\alpha \stackrel{w^*} \longrightarrow c $ in $\mathfrak{A} $. Because $c_\alpha a  \stackrel{wk} \longrightarrow c a $, for every $\varphi \in \mathfrak{A}^*$ we have $$ lim_\alpha \langle \varphi , c_\alpha (a b ) \rangle      =   lim_\alpha    \langle  b \cdot \varphi , c_\alpha a \rangle  =  \langle  b \cdot \varphi , c a \rangle =  \langle \varphi , c (a b ) \rangle $$ so that $ a b \in \sigma wc (\mathfrak{A} )$. Next, by $w^*$-continuity of the multiplication, $ c_\alpha  b \stackrel{w^*} \longrightarrow c b$. Then $c_\alpha ( b a) = (c_\alpha b ) a  \stackrel{wk} \longrightarrow ( c b ) a = c  (  b a)$, which means $ b a \in \sigma wc (\mathfrak{A} )$.
\end{rem}

For a given dual Banach algebra $\mathfrak{A} $, it would be interesting to determine the set  $\sigma wc (\mathfrak{A} )$, see for instance Proposition \ref{9.1} and Example \ref{9.2} below. However, special care should be taken with the trivial cases $\sigma wc (\mathfrak{A} ) = \{ 0 \}$ and $\sigma wc (\mathfrak{A} ) = \mathfrak{A}$, as follows.

\begin{rem} \label{7.2} Let $\mathfrak{A}$ be a dual Banach algebra. Then: 

$(i)$ If  $\sigma wc (\mathfrak{A} ) = \{ 0 \}$, then $j_{\mathfrak{A}} $ is the zero map. Therefore, for every derivation $D
:\mathfrak{A} \longrightarrow \mathfrak{A}^*$, $j_{\mathfrak{A}} \circ D = 0$. So, $\mathfrak{A} $ is weakly Connes amenable.

$(ii)$ If $\sigma wc (\mathfrak{A} ) = \mathfrak{A}$ (or equivalently $\mathfrak{A}^*$ is normal by \cite[Proposition 4.4]{r3}), then  $j_{\mathfrak{A}} $ is the identity map on $ \mathfrak{A}^*$ and then $j_{\mathfrak{A}} \circ D = D$. In particular, Definition \ref{2.1} becomes: $\mathfrak{A} $ is weakly Connes amenable if and only if every $w^*$-$w^*$ continuous derivation $D : \mathfrak{A} \longrightarrow \mathfrak{A}^*$ is inner.
\end{rem}

As the name suggests, weak Connes amenability is weaker than Connes amenability.
\begin{theo} \label{0.2} Every Connes amenable dual Banach algebra
is weakly Connes amenable.
\end{theo}
{\bf Proof.} Let $\mathfrak{A} = (\mathfrak{A}_*)^*$ be a Connes amenable dual Banach algebra, and let  $D :\mathfrak{A} \longrightarrow \mathfrak{A}^*$ be a derivation such that $j_{\mathfrak{A}} \circ D : \mathfrak{A} \longrightarrow   \sigma wc (\mathfrak{A} )^*$ is $w^*$-$w^*$ continuous. By the assumption, $j_{\mathfrak{A}} \circ D$ is inner, as required. \qed

For dual Banach algebras, weak amenability implies weak Connes amenability as follows.

\begin{theo} \label{2.12} Every weakly amenable dual Banach algebra
is weakly Connes amenable.
\end{theo}
{\bf Proof.} Let $\mathfrak{A} = (\mathfrak{A}_*)^*$ be a weakly amenable dual Banach algebra, and let  $D :\mathfrak{A} \longrightarrow \mathfrak{A}^*$ be a derivation such that $j_{\mathfrak{A}} \circ D : \mathfrak{A} \longrightarrow   \sigma wc (\mathfrak{A} )^*$ is $w^*$-$w^*$ continuous. By the assumption, $D$ is inner. Then by Lemma \ref{0.1}$(ii)$,  $j_{\mathfrak{A}} \circ D$ is inner. \qed

Suppose that $\mathfrak{A}$ is a Banach algebra and that $E$ is a
Banach $\mathfrak{A}$-bimodule. A derivation $D:\mathfrak{A} \longrightarrow E^*$ is $w^*$-\textit{approximately inner} if there exists a net $(f_\alpha)_\alpha \subseteq E^*$ such that $ D (a) =w^*-\lim_\alpha ad_{f_\alpha} ( a )$ for all $a \in
\mathfrak{A}$ \cite{GL}. The concept of $w^*$-approximate Connes amenability were introduced in \cite{m}. A dual Banach algebra $\mathfrak{A}$ is $w^*$-\textit{approximately Connes amenable}
if for every normal, dual Banach $\mathfrak{A}$-bimodule $E$, every
$w^*$-$w^*$ continuous derivation $D : \mathfrak{A} \longrightarrow
E$ is $w^*$-approximately inner.

\begin{defi} \label{2.5} A dual Banach algebra $\mathfrak{A} = (\mathfrak{A}_*)^*$ is $w^*$-\textit{approximately weakly Connes amenable} if for every derivation $D :\mathfrak{A} \longrightarrow \mathfrak{A}^*$ such that  $j_{\mathfrak{A}} \circ D : \mathfrak{A} \longrightarrow   \sigma wc (\mathfrak{A} )^*$ is $w^*$-$w^*$ continuous, derivation  $j_{\mathfrak{A}} \circ D$ is $w^*$-approximately inner.
\end{defi}

\begin{prop} \label{0.4} Every $w^*$-approximately Connes amenable dual Banach
algebra is $w^*$-approximately weakly Connes amenable.
\end{prop}
{\bf Proof.} The proof is analogous to that of Theorem \ref{0.2}. \qed

Let $\mathfrak{A}$ be a dual Banach algebra. Composing the canonical inclusion $ \mathfrak{A} \hat{\otimes} \mathfrak{A}
\longrightarrow ( \mathfrak{A} \hat{\otimes} \mathfrak{A})^{**}$ with the quotient map $ ( \mathfrak{A} \hat{\otimes} \mathfrak{A})^{**}
\longrightarrow \sigma wc ( ( \mathfrak{A} \hat{\otimes} \mathfrak{A})^*)^*$, we obtain an
$\mathfrak{A}$-bimodule homomorphism $\zeta : \mathfrak{A} \hat{\otimes} \mathfrak{A} \longrightarrow
\sigma wc ( (\mathfrak{A} \hat{\otimes} \mathfrak{A}) ^*)^*$ with $w^*$-dense range. 

\begin{defi} \label{9.3} (\cite{m}) A dual Banach algebra $\mathfrak{A}$
is \textit{pseudo Connes amenable} if there exists a net
$(m_\alpha)_\alpha$ in $ \mathfrak{A} \hat{\otimes} \mathfrak{A}$, called an \textit{approximate} $\sigma
wc$-\textit{diagonal} for $\mathfrak{A}$, such that
 $ a \ . \  \zeta (m_\alpha) - \zeta (m_\alpha) \ . \
a \stackrel{w^*} \longrightarrow 0 $, and $a
\pi_{\sigma wc} (m_\alpha) \stackrel{w^*} \longrightarrow a$ for every $a \in \mathfrak{A}$.
\end{defi}

We need the following, which is \cite[Proposition 5.6]{m}.
\begin{prop} \label{0.5} Let $\mathfrak{A}$ be a pseudo Connes amenable dual
Banach algebra, and let $E$ be a normal, dual Banach $\mathfrak{A}$-bimodule such that each $w^*$-approximate identity of $\mathfrak{A}$ is also a one-sided $w^*$-approximate identity for $E$. Then every $w^*$-$w^*$ continuous derivation $D :\mathfrak{A}  \longrightarrow E$ is $w^*$-approximately inner.
\end{prop}

We now discuss relations between pseudo Connes amenability and these
new notions.

\begin{theo} \label{2.7} Every pseudo Connes amenable dual
Banach algebra is $w^*$-approximately weakly  Connes
amenable.
\end{theo}
{\bf Proof.} Let $\mathfrak{A} = (\mathfrak{A}_*)^*$ be a pseudo Connes amenable dual Banach algebra, and let  $D :\mathfrak{A} \longrightarrow \mathfrak{A}^*$ be a derivation such that $j_{\mathfrak{A}} \circ D : \mathfrak{A} \longrightarrow   \sigma wc (\mathfrak{A} )^*$ is $w^*$-$w^*$ continuous. We use Proposition \ref{0.5} with $\sigma wc (\mathfrak{A} )^*$ in place of $E$. Thus $j_{\mathfrak{A}} \circ D$  is $w^*$-approximately inner, as required. \qed

\begin{cor} \label{2.8}  Every commutative, pseudo Connes amenable dual
Banach algebra is weakly Connes amenable.
\end{cor}
{\bf Proof.} Let $\mathfrak{A} = (\mathfrak{A}_*)^*$ be a pseudo Connes amenable dual Banach algebra, and let  $D :\mathfrak{A} \longrightarrow \mathfrak{A}^*$ be a derivation such that $j_{\mathfrak{A}} \circ D : \mathfrak{A} \longrightarrow   \sigma wc (\mathfrak{A} )^*$ is $w^*$-$w^*$ continuous. By Theorem \ref{2.7}, there is a net $ (\varphi_\alpha)_\alpha$ in $\sigma wc (\mathfrak{A} )^*$ that $j_{\mathfrak{A}} \circ D (a) =w^*-\lim_\alpha ad_{\varphi_\alpha}(a)$ for $a \in \mathfrak{A}$. As $\mathfrak{A}$ is commutative, $ ad_{\varphi} =0$ for each $ \varphi \in \sigma wc (\mathfrak{A} )^*$. Hence $j_{\mathfrak{A}} \circ D = 0$, and then $\mathfrak{A}$ is weakly  Connes amenable. \qed

\begin{prop} \label{6.1} Every commutative, $w^*$-approximately Connes amenable dual
Banach algebra is weakly  Connes amenable.
\end{prop}
{\bf Proof.} The proof is similar to that of Corollary \ref{2.8}. \qed

Suppose that $\mathfrak{A}=(\mathfrak{A}_*)^*$ and $\mathfrak{B} =
(\mathfrak{B}_*)^*$ are dual Banach algebras. We consider the
$\ell^1$-direct sum $\mathfrak{A} \oplus^1 \mathfrak{B}$ with norm $
|| (a , b) || = || a|| + ||b||$ for $a \in \mathfrak{A}$ and $b \in
\mathfrak{B}$. This is a dual Banach algebra under pointwise-defined
operations and with predual the $\ell^\infty$-direct sum
$\mathfrak{A}_* \oplus^\infty \mathfrak{B}_*$, where the norm
$||.||_\infty$ is defined through $ ||(\phi , \psi )||_\infty = \max
( ||\phi|| , ||\psi||)$ for $\phi \in \mathfrak{A}_*$ and $\psi \in
\mathfrak{B}_*$. The duality is given by $$ \langle  (a , b) , (\phi , \psi )  \rangle = \langle   a , \phi \rangle + \langle
 b , \psi  \rangle \ \ \ \ ( a \in \mathfrak{A}, \ b \in \mathfrak{B},
\ \phi \in \mathfrak{A}_*, \ \psi \in \mathfrak{B}_*) \ .$$ We write
$ \imath_{\mathfrak{A}} : \mathfrak{A} \longrightarrow \mathfrak{A} \oplus^1 \mathfrak{B}$ for the natural injective homomorphism. 
 It is obvious that
$ \imath_{\mathfrak{A}}^*( \varphi , \psi) = \varphi$ for $ (
\varphi , \psi) \in \mathfrak{A}^* \oplus^\infty \mathfrak{B}^* =
(\mathfrak{A} \oplus^1 \mathfrak{B} )^*$.  It follows from \cite[Lemma 2.6]{m} that $  \sigma wc (\mathfrak{A} \oplus \mathfrak{B} )= \sigma wc (\mathfrak{A})\oplus  \sigma w c (\mathfrak{B})$. We also denote by $ \nu_{\mathfrak{A}} : \sigma wc (\mathfrak{A} \oplus^1 \mathfrak{B} )^* =  \sigma wc (\mathfrak{A} )^*  \oplus^\infty   \sigma wc (\mathfrak{B} )^* \longrightarrow  \sigma wc (\mathfrak{A} )^*$ the adjoint of the natural embedding $ \sigma wc (\mathfrak{A} ) \hookrightarrow  \sigma wc (\mathfrak{A} ) \oplus^1 \sigma wc (\mathfrak{B} ) =  \sigma wc (\mathfrak{A} \oplus^1 \mathfrak{B} )$, so that $ \nu_{\mathfrak{A}} ( \varphi , \psi) = \varphi$ for $ (
\varphi , \psi) \in  \sigma wc (\mathfrak{A} ) ^* \oplus^\infty  \sigma wc (\mathfrak{B} ) ^* $. Similarly we define $
\imath_{\mathfrak{B}}$ and $\nu_{\mathfrak{B}}$.

\begin{lemm} \label{0.3} Let $\mathfrak{A}$ and $\mathfrak{B} $ be  dual Banach algebras, and let $D :\mathfrak{A} \oplus^1 \mathfrak{B} \longrightarrow (\mathfrak{A} \oplus^1 \mathfrak{B} )^* = \mathfrak{A}^* \oplus^\infty \mathfrak{B}^*$ be a derivation. Then:

$(i)$ $j_{\mathfrak{A}} \circ ( \imath_{\mathfrak{A}}^* \circ D  ) =\nu_{\mathfrak{A}} \circ (j_{\mathfrak{A} \oplus^1 \mathfrak{B}} \circ D ) $;

$(ii)$  $j_{\mathfrak{B}} \circ ( \imath_{\mathfrak{B}}^* \circ D  ) =\nu_{\mathfrak{B}} \circ (j_{\mathfrak{A} \oplus^1 \mathfrak{B}} \circ D ) $.
\end{lemm}
{\bf Proof.} We only prove $(i)$. For every $a \in \mathfrak{A}$, $b \in \mathfrak{B} $, and $c \in \sigma wc (\mathfrak{A}) $ we have
\begin{align*}
\langle  c , j_{\mathfrak{A}} \circ ( \imath_{\mathfrak{A}}^* \circ D  ) (a, b)  \rangle &= \langle  c ,  ( \imath_{\mathfrak{A}}^* \circ D  ) (a, b)   \rangle = \langle (c , 0 ) , D (a , b)  \rangle \\&=  \langle   (c , 0 ) , (j_{\mathfrak{A} \oplus^1 \mathfrak{B}} \circ D   ) (a , b)  \rangle =
 \langle  c  , \nu_{\mathfrak{A}} \circ (j_{\mathfrak{A} \oplus^1 \mathfrak{B}} \circ D   ) (a , b)  \rangle \ ,
\end{align*}
as required. \qed

\begin{theo} \label{5.2} Suppose that $\mathfrak{A}$ and $\mathfrak{B} $ are weakly Connes amenable dual
Banach algebras such that $\mathfrak{A}^2$ and $\mathfrak{B}^2$ are $w^*$-dense in $\mathfrak{A}$ and $\mathfrak{B}$, respectively. Then $\mathfrak{A} \oplus^1 \mathfrak{B}$ is weakly Connes amenable.
\end{theo}
{\bf Proof.} Take a derivation $D :\mathfrak{A} \oplus^1 \mathfrak{B} \longrightarrow (\mathfrak{A}
\oplus^1 \mathfrak{B} )^* = \mathfrak{A}^* \oplus^\infty \mathfrak{B}^*$ for which $j_{\mathfrak{A} \oplus^1 \mathfrak{B}} \circ D :\mathfrak{A} \oplus^1 \mathfrak{B} \longrightarrow   \sigma wc (\mathfrak{A} \oplus^1 \mathfrak{B} )^* =  \sigma wc (\mathfrak{A} )^*  \oplus^\infty   \sigma wc (\mathfrak{B} )^* $ is $w^*$-$w^*$ continuous. Then  $ \imath_{\mathfrak{A}}^* \circ D \circ \imath_{\mathfrak{A}} :\mathfrak{A} \longrightarrow \mathfrak{A}^*$
is a derivation. An easy calculation shows that $j_{\mathfrak{A}} \circ ( \imath_{\mathfrak{A}}^* \circ D \circ \imath_{\mathfrak{A}} ) : \mathfrak{A} \longrightarrow   \sigma wc (\mathfrak{A} )^*$ is $w^*$-$w^*$ continuous. Since $\mathfrak{A}$ is weakly Connes amenable, there exists $ \varphi \in \sigma wc (\mathfrak{A} )^*$ such that $j_{\mathfrak{A}} \circ ( \imath_{\mathfrak{A}}^* \circ D)
(a , 0) = ad_\varphi (a)$ and then, by Lemma \ref{0.3}$(i)$, $\nu_{\mathfrak{A}} \circ (j_{\mathfrak{A} \oplus^1 \mathfrak{B}} \circ D ) (a , 0) = ad_\varphi  (a)$ for every $a \in \mathfrak{A}$.  We regard $\varphi = (\varphi , 0 )$ as an element of $\sigma wc (\mathfrak{A} )^* \oplus^\infty \sigma wc (\mathfrak{B} )^* $. Then it is easily checked that $\nu_{\mathfrak{A}} \circ (j_{\mathfrak{A} \oplus^1 \mathfrak{B}} \circ ( D - ad_\varphi ) ) (a , 0) = 0$. So, by
replacing $D$ by $D-ad_\varphi$, we may suppose that  $\nu_{\mathfrak{A}} \circ (j_{\mathfrak{A} \oplus^1 \mathfrak{B}} \circ D ) (a , 0) =0$  for every $a \in \mathfrak{A}$. For each  $ a_1 , a_2 \in \mathfrak{A}$, $ a \in \sigma wc (\mathfrak{A} )$ and $b \in  \sigma wc (\mathfrak{B} )$ we see that
\begin{align*}
\langle (a , b) , (j_{\mathfrak{A} \oplus^1 \mathfrak{B}}  \circ D  ) (a_1 a_2, 0)  \rangle &= \langle (a , b) , j_{\mathfrak{A} \oplus^1 \mathfrak{B}} ( (  a_1 , 0) \cdot D(a_2 , 0 ) +  D(a_1 , 0 ) \cdot (a_2 , 0 )  ) \rangle \\&= \langle (a , b) ,  (  a_1 , 0) \cdot (j_{\mathfrak{A} \oplus^1 \mathfrak{B}} \circ D )  (a_2 , 0 ) + (j_{\mathfrak{A} \oplus^1 \mathfrak{B}} \circ D )  (a_1 , 0 ) \cdot  (a_2 , 0 ) \rangle \\&= \langle (a   a_1 , 0),   (j_{\mathfrak{A} \oplus^1 \mathfrak{B}} \circ D )  (a_2 , 0 ) + ( a_2 a  , 0),   (j_{\mathfrak{A} \oplus^1 \mathfrak{B}} \circ D )  (a_1 , 0 ) \rangle \\&=  \langle a   a_1,   \nu_{\mathfrak{A}} \circ (j_{\mathfrak{A} \oplus^1 \mathfrak{B}} \circ D )  (a_2 , 0 ) \rangle + \langle  a_2  a,   \nu_{\mathfrak{A}} \circ (j_{\mathfrak{A} \oplus^1 \mathfrak{B}} \circ D )  (a_1 , 0 ) \rangle \\&= 0 + 0 = 0 , \ \  ( a   a_1,  a_2  a \in \sigma wc (\mathfrak{A} ) \  \text{by Remark} \ \ref{7.1} )
\end{align*}
and whence $  (j_{\mathfrak{A} \oplus^1 \mathfrak{B}}  \circ D  ) (a_1 a_2, 0)  = 0$ for all $ a_1 , a_2 \in \mathfrak{A}$. Now, it follows from $w^*$-density of $\mathfrak{A}^2$ and $w^*$-continuity of $ j_{\mathfrak{A} \oplus^1 \mathfrak{B}}  \circ D $ that $  (j_{\mathfrak{A} \oplus^1 \mathfrak{B}}  \circ D  ) (a, 0)  = 0$ for every $ a \in \mathfrak{A}$. 

A similar argument, with $\mathfrak{B} $ in place of $\mathfrak{A} $, shows that $  (j_{\mathfrak{A} \oplus^1 \mathfrak{B}}  \circ D  ) (0, b)  = 0$ for every $ b \in \mathfrak{B}$. It follows that $ j_{\mathfrak{A} \oplus^1 \mathfrak{B}}  \circ D  = 0$, and hence $\mathfrak{A} \oplus^1 \mathfrak{B}$ is weakly Connes amenable. \qed

\begin{prop} \label{2.3} Let $\mathfrak{A}$ be a commutative Banach algebra, let $\mathfrak{B} $ be a
dual Banach algebra, and let $\theta: \mathfrak{A} \longrightarrow
\mathfrak{B}$ be a (continuous) homomorphism with $w^*$-dense range.
If $\mathfrak{A}$ is weakly amenable, then
$\mathfrak{B} $ is weakly Connes amenable.
\end{prop}

{\bf Proof.} Take a derivation $D
:\mathfrak{B} \longrightarrow \mathfrak{B}^* $ such that $j_{\mathfrak{B}}  \circ D : \mathfrak{B} \longrightarrow
   \sigma wc (\mathfrak{B} )^* $ is $w^*$-$w^*$ continuous. Then $(j_{\mathfrak{B}}  \circ D ) \circ \theta :
\mathfrak{A} \longrightarrow  \sigma wc (\mathfrak{B} )^*$ is a derivation. By
\cite[Theorem 2.8.63 (iii)]{dal}, $(j_{\mathfrak{B}}  \circ D ) \circ \theta= 0$. From $w^*$-continuity of $j_{\mathfrak{B}}  \circ D$ and $w^*$-density of the range of $\theta$, we conclude that  $j_{\mathfrak{B}}  \circ D = 0$, as required. \qed

Let $\mathfrak{A}$ be a Banach algebra. We write $ WAP
(\mathfrak{A}^*)$ for the space of all \textit{weakly almost
periodic functionals} on $\mathfrak{A}$. It is noted by V. Runde
that $ WAP (\mathfrak{A}^*)^*$ is a dual Banach algebra with a
universal property \cite[Theorem 4.10]{r3}. Later, M. Daws called $
WAP (\mathfrak{A}^*)^*$ the \textit{dual Banach algebra enveloping
algebra of} $\mathfrak{A}$ \cite[Definition 2.10]{daw2}. There is a (continuous) homomorphism $ \kappa : \mathfrak{A} \longrightarrow  WAP (\mathfrak{A}^*)^*$ whose range is  $w^*$-dense. Indeed, the map $ \kappa$ obtained by composing the canonical inclusion $ \mathfrak{A} \longrightarrow \mathfrak{A}^{**}$ with the adjoint of the inclusion map $ WAP (\mathfrak{A}^*) \hookrightarrow  \mathfrak{A} ^*$  \cite{r3}.

\begin{cor} \label{2.11} Let $\mathfrak{A}$ be a commutative, weakly
amenable Banach algebra. Then $ WAP (\mathfrak{A}^*)^*$ is weakly Connes amenable.
\end{cor}
{\bf Proof.} We apply Proposition \ref{2.3} for the canonical homomorphism $ \kappa : \mathfrak{A} \longrightarrow  WAP (\mathfrak{A}^*)^*$. \qed

\begin{cor} \label{2.14} Let $\mathfrak{A}$ be a commutative $C^*$-algebra. Then $ WAP (\mathfrak{A}^*)^*$ is weakly Connes  amenable.
\end{cor}
{\bf Proof.} Since $C^*$-algebras are weakly amenable \cite{h}, it
is immediate by Proposition \ref{2.3}. \qed

\section{Examples}

In this section we give some examples to show difference between
weak Connes amenability and some older notions such as Connes amenability,
weak amenability and pseudo (Connes) amenability.

\begin{ex} \label{2.10}  It is well known that $\ell^1= \ell^1( \mathbb{N})$ with pointwise multiplication is a commutative,
weakly amenable dual Banach algebra. It follows from Proposition
\ref{2.12} that $\ell^1$ is weakly Connes amenable. Further,
$WAP((\ell^1)^*)^* = WAP(\ell^\infty)^*$ is also weakly Connes 
amenable, by Corollary \ref{2.11}. However, they fail to be Connes
amenable because of the lack of an identity.
\end{ex}

\begin{ex}  \label{2.9} Consider the semigroup $\mathbb{N}_{\vee}$
 which is $\mathbb{N }$ with the operation
$m \vee n := \max \{m, n\}$,  $(m, n \in \mathbb{N})  $. It is known
that $\ell^1( \mathbb{N}_{\vee})$ is a commutative, weakly amenable
dual Banach algebra. An argument similar to Example \ref{2.10} shows
that both $\ell^1( \mathbb{N}_{\vee})$ and $WAP(\ell^1(
\mathbb{N}_{\vee})^*)^*$ are weakly Connes amenable. Notice that
none of them are Connes amenable by \cite[Theorem 5.13]{daw1} and
\cite[Page 262]{daw2}, respectively.
\end{ex}
\begin{ex} \label{13.2} We denote by $\mathbb{N}_{\wedge}$
the semigroup $\mathbb{N }$ with the operation $m \wedge n := \min
\{m, n\}$,  $(m, n \in \mathbb{N})  $. It is known that $\ell^1(
\mathbb{N}_{\wedge})$ is a commutative, weakly amenable Banach
algebra. By Corollary \ref{2.11}, $WAP(\ell^1( \mathbb{N}_{\wedge})^*)^*$ is weakly Connes amenable. It was shown in \cite[Theorem 7.6]{daw2} that $WAP(\ell^1( \mathbb{N}_{\wedge})^*)^*$ is not Connes amenable.
\end{ex}

It is well known that the measure algebra $M(G)$ of a locally compact
group $G$, is a dual Banach algebra with predual $C_0(G)$. It was shown that $M(G)$ is Connes amenable if and only if $G$ is amenable \cite{r4}. $M(G)$ is weakly amenable if and only if $G$ is discrete, see for
instance \cite[Theorem 4.2.13]{r2}.
\begin{ex} \label{8.1} Let $G$ be a locally compact group. It follows from Theorems \ref{2.12} and \ref{0.2} that $M(G)$ is weakly Connes amenable if $G$ is either discrete or amenable. It would be desirable to characterize those locally compact
groups $G$ for which $M(G)$ is weakly Connes amenable.
\end{ex}

\begin{theo} \label{5.4} Let $G$ be a non-discrete abelian locally compact group, and set $ \mathfrak{A} := M(G) \oplus^1
\ell^1$. Then:

$(i)$ $ \mathfrak{A} $ is not weakly amenable;

$(ii)$ $ \mathfrak{A} $ is not Connes amenable;

$(iii)$ $ \mathfrak{A} $ is weakly Connes amenable.
\end{theo}
{\bf Proof.} $(i)$ If $ \mathfrak{A}$ were weakly amenable, then so is its image under the natural projection 
$ \mathcal{P}_{M(G)} : \mathfrak{A} \longrightarrow  M(G) $ by \cite[Proposition
2.8.64]{dal}. However, this is not possible for non-discrete groups \cite[Theorem 4.2.13]{r2}.

$(ii)$  If $ \mathfrak{A}$ were Connes amenable, then so is its image under the natural projection 
$ \mathcal{P}_{\ell^1} : \mathfrak{A} \longrightarrow  \ell^1 $ by \cite[Proposition 4.2]{r1}. But, as mentioned in Example \ref{2.10}, $\ell^1$ is not Connes amenable.

$(iii)$ Both $\ell^1 $ and $M(G) $ are weakly Connes amenable by Examples \ref{2.10} and \ref{8.1}, respectively. The result now follows from Theorem \ref{5.2}. \qed

\begin{rem} \label{5.5} Suppose that $G$ is either a discrete or an amenable locally compact group. The same argument as in Theorem \ref{5.4} $(iii)$ shows that  $M(G) \oplus^1 \ell^1$ is weakly Connes amenable.
\end{rem}

We write $ \mathfrak{A}^\sharp$ for the unitization of an algebra $
\mathfrak{A}$.
\begin{theo} \label{5.3} Let $G$ be a non-discrete abelian locally compact group, and set $ \mathfrak{A} := M(G) \oplus^1
(\ell^1)^\sharp$. Then:

$(i)$ $ \mathfrak{A} $ is not weakly amenable;

$(ii)$ $ \mathfrak{A} $ is not pseudo amenable;

$(iii)$ $ \mathfrak{A} $ is not Connes amenable;

$(iv)$ $ \mathfrak{A} $ is not pseudo Connes amenable;

$(v)$ $ \mathfrak{A} $ is weakly Connes amenable.
\end{theo}
{\bf Proof.} $(i)$ This is essentially the same as the proof of
Theorem \ref{5.4} $(i)$.

 $(ii)$ By \cite[Theorem 3.1]{GZ} and
\cite[Theorem 4.1]{dlz}, $(\ell^1)^\sharp$ is not pseudo amenable.
As $(\ell^1)^\sharp$ is the image of $ \mathfrak{A}$ under the projection $ \mathcal{P}_{(\ell^1)^\sharp} : \mathfrak{A} \longrightarrow  (\ell^1)^\sharp $, $ \mathfrak{A}$ can not be pseudo
amenable by \cite[Proposition 2.2]{GZ}.

$(iii)$  Since $(\ell^1)^\sharp$ is not Connes amenable, a similar argument as in Theorem \ref{5.4} $(ii)$ holds.

 $(iv)$ It follows from \cite[Theorem 3.2 and Theorem 5.1]{m} that $(\ell^1)^\sharp$ is
 not pseudo Connes amenable. Towards a contradiction, suppose that  $ \mathfrak{A}$ is pseudo Connes
 amenable. Then $(\ell^1)^\sharp$, the image of $ \mathfrak{A}$
 under $\mathcal{P}_{(\ell^1)^\sharp}$, must be pseudo Connes amenable by
\cite[Proposition 4.5]{m}, which is not the case.

$(v)$ Weak amenability of $\ell^1$ implies that of
$(\ell^1)^\sharp$. Using Proposition \ref{2.12} and Example \ref{8.1}, weak Connes amenability of $ \mathfrak{A}$ is a
consequence of Theorem \ref{5.2}. \qed

Let $\mathcal{V}$ be a Banach space of dimension (at least) 2, and let $f \in \mathcal{V}^*$
be a non-zero element such that $|| f || \leq 1$. Then $\mathcal{V}$
equipped with the product defined by $ a b := f(a) b$ for $a , b
\in\mathcal{V}$, is a non-commutative Banach algebra denoted by $\mathcal{V}_{f} $.  Moreover, if $\mathcal{V}$ is a dual Banach space and if $f$ is $w^*$-continuous, then $\mathcal{V}_{f}$ is a dual Banach
algebra, see for instance \cite{zsm}. 

\begin{prop} \label{9.1}Let $\mathcal{V}$ be a dual Banach space with $ dim \mathcal{V} \geq 2 $, and let $0 \neq f \in \mathcal{V}^*$
be $w^*$-continuous such that $|| f || \leq 1$. Then:

$(i)$  $\sigma wc (\mathcal{V}_{f} ) =  ker f$;

$(ii)$ $\mathcal{V}_{f} $ is not pseudo Connes amenable;

$(iii)$  $\mathcal{V}_{f} $ is not $w^*$-approximately Connes amenable;

$(iv)$ $\mathcal{V}_{f}$ is weakly Connes amenable.
\end{prop}
{\bf Proof.} $(i)$ It is easy to verify that every left ideal of $\mathcal{V}_{f}$ is contained in $ker f$. In particular, $\sigma wc (\mathcal{V}_{f} ) \subseteq ker f$. For the converse, first note that the map $\mathcal{V}_{f}  \longrightarrow  \mathcal{V}_{f} $, $ a \longmapsto ab =  f(a) b$ is $w^*$-$wk$ continuous for all $b \in \mathcal{V}_{f}$. Now, take an arbitrary element $ b \in ker f $ and let $ a_\alpha \stackrel{w^*} \longrightarrow a$ in $\mathcal{V}_{f}$. Then, for all $ \varphi \in \mathcal{V}_{f}^*$ and all $\alpha$ we have $  \langle \varphi , b a_\alpha \rangle =  f(b)  \langle \varphi , a_\alpha \rangle = 0 =  \langle \varphi ,  b a \rangle $, i.e., $b a_\alpha  \stackrel{wk} \longrightarrow b a$. Thus $ b \in \sigma wc (\mathcal{V}_{f} )$, as required. 

$(ii)$ Assume that $\mathcal{V}_{f} $ is pseudo Connes
amenable and that $(m_\alpha) \subseteq \mathfrak{A} \hat{\otimes}
\mathfrak{A}$ is an approximate $\sigma wc$-diagonal for
$\mathfrak{A}$. Take a non-zero element $a \in \ker f$. Then $ 0 = f(a)
\pi_{\sigma wc} (m_\alpha)  = a \pi_{\sigma wc} (m_\alpha)
\stackrel{w^*} \longrightarrow a$, which implies that $ \langle x, a
\rangle = 0$ for each $x$ in the predual of $\mathcal{V}_{f} $. It
forces $a$ to be zero, a contradiction.

$(iii)$ Towards a contradiction, suppose that $\mathcal{V}_{f} $ is $w^*$-approximately
Connes amenable. By \cite[Proposition 2.2]{m}, $\mathcal{V}_{f} $
has a right $w^*$-approximate identity, i.e., there is a net
$(e_\alpha) \subseteq \mathcal{V}_{f} $ for which $f(a) e_\alpha = a
e_\alpha \stackrel{w^*} \longrightarrow a$ for every $a \in
\mathcal{V}_{f}$. We take a non-zero element $a \in \ker f$, and then the rest of the proof is analogous to $(ii)$.

$(iv)$  It was shown that $\mathcal{V}_{f}$ is weakly amenable \cite[Page 507]{z}, see also \cite[Proposition 2.13]{dls}. Hence, $\mathcal{V}_{f}$ is automatically weakly Connes amenable. \qed

\begin{ex} \label{9.2}It is denoted by $ \mathbb{M}_n$ the collection of all $ n \times n$ matrices, $n \in \mathbb{N}$, 
with entries from $\mathbb{C}$. It is known that $ \mathbb{M}_n =  \mathbb{M}_n^*$ is a dual Banach
algebra, and that $w^*$-topology and $wk$-topology are the same on $ \mathbb{M}_n$. It is then straightforward that $ \sigma wc (\mathbb{M}_n)= \mathbb{M}_n$.
\end{ex}


\begin{thebibliography}{99}
\bibitem{bcd}
W. G. Bade, P. C. Curtis, H. G. Dales, Amenability and weak
amenability for Beurling and Lipschitz algebras. {\it Proc. London
Math. Soc.} (3) {\bf 55} (1987), 359-377.
\bibitem{dal}
H. G. Dales, Banach algebras and automatic continuity, {\it
Clarendon Press}, Oxford, 2000.


\bibitem{dls}
 H. G. Dales, A. T.-M. Lau, D. Strauss, Banach algebras on semigroups
and on their compactifications. {\it Mem. Amer. Math. Soc.} {\bf 205} (966) (2010).

\bibitem{dlz}
 H. G. Dales, R. J. Loy, Y. Zang, Approximate amenability for
Banach sequence algebras. {\it Studia Math.} {\bf 177} (2006),
81-96.
\bibitem{daw1}
 M. Daws, Connes-amenability of bidual and weighted semigroup
algebras, {\it Math. Scand.} {\bf 99} (2006), 217-246.
\bibitem{daw2}
 M. Daws, Dual Banach algebras: representations and injectivity,
{\it Studia Math.} {\bf 178} (2007), 231-275.



\bibitem{GL}
 F. Ghahramani, R. J. Loy, Generalized notions of amenability. {\it J. Func. Anal.} {\bf 208} (2004), 229-260.

\bibitem{GZ}
 F. Ghahramani, Y. Zhang, Pseudo-amenable and
pseudo-contractible Banach algebras. {\it Math. Proc. Camb. Phil.
Soc.} {\bf 142} (2007), 111-123.
\bibitem{h}
U. Haagerup, All nuclear $C^*$-algebras are amenable. {\it Invent.
Math.} {\bf 74} (1983), 305-319.


\bibitem{joh1}
 B. E. Johnson, Derivations from $ L^1(G)$ into $ L^1(G)$ and $
 L^\infty(G)$. In: J. P. Pier (ed.), \textit{Harmonic Analysis(Luxembourg,
 1987)}, pp. 191-198. \textit{Springer Verlag}, 1988.


\bibitem{joh}
 B. E. Johnson, Cohomology in Banach algebras, {\it Mem. Amer.
Math. Soc.} {\bf 127} (1972).

\bibitem{m}
A. Mahmoodi, Connes-amenability-like properties, {\it Studia Math.}
{\bf 220} (2014), 55-72.
\bibitem{r1}
 V. Runde,  Amenability for dual Banach algebras, {\it Studia
Math.} {\bf 148} (2001), 47-66.

\bibitem{r4}
 V. Runde, Connes-amenability and normal,
virtual diagonals for measure algebras I. {\it J. London Math. Soc.}
{\bf 67} (2003), 643-656.
\bibitem{r3}
 V. Runde, Dual Banach algebras: Connes-amenability, normal,
virtual diagonals, and injectivity of the predual bimodule, {\it
Math. Scand.} {\bf 95} (2004), 124-144.

\bibitem{r2}
 V. Runde, Amenable Banach Algebras: A Panorama, {\it Springer Monographs in Mathematics}, 2020.

\bibitem{z}
Y. Zhang, Weak amenability of a class of Banach algebras, {\it Canad. Math. Bull.} {\bf 44} (2001), 504-508.

\bibitem{zsm}
 M. Ziamanesh, B. Shojaee, A. Mahmoodi, Essential amenability of dual Banach algebras, {\it
Bull. Aust. Math. Soc.} {\bf 100} (2019), 479-488.


\end{thebibliography}
\end{document}